\documentclass
{article}
\usepackage[T2A]{fontenc}
\usepackage[cp1251]{inputenc}
\usepackage[english,russian]{babel}
\usepackage[tbtags]{amsmath}
\usepackage{amsfonts,amssymb,mathrsfs,amscd}
\usepackage{wasysym}
\usepackage{verbatim}
\usepackage{eucal}
\usepackage{upgreek}
\usepackage{graphicx}
\usepackage{enumerate}
\let\No=\textnumero

\usepackage[hyper]{izv2e}
\JournalName{}
\numberwithin{equation}{section}

\theoremstyle{plain}
\newtheorem{theorem}{Теорема}
\newtheorem{maintheorem}{Основная теорема}

\newtheorem{ThSIw}{Теорема о малых интервалах с весом}

\theoremstyle{definition}
\newtheorem{definition}{Определение}

\newtheorem{example}{Пример}

\renewcommand{\leq}{\leqslant} 
\renewcommand{\geq}{\geqslant}
\newcommand{\RR}{\mathbb{R}} 
\newcommand{\CC}{\mathbb{C}}

\newcommand{\rad}{\text{\tiny\rm rd}}

\DeclareMathOperator{\dd}{\,{\mathrm  d\!}}

\DeclareMathOperator{\mes}{mes}

\begin{document} 
\title{Мероморфные функции и разности субгармонических функций в интегралах  и разностная характеристика Неванлинны. II. Явные оценки интеграла от радиальной максимальной характеристики роста}


\author[B.\,N.~Khabibullin]{Б.\,Н.~Хабибуллин}
\address{Башкирский государственный университет}
\email{khabib-bulat@mail.ru}

\date{}
\udk{517.547.26 : 517.547.28 : 517.574}

 \maketitle

\begin{fulltext}

\begin{abstract} Пусть  $U\not\equiv \pm\infty$ --- разность субгармонических функций, т.е. $\delta$-субгар\-м\-о\-н\-и\-ч\-е\-ская  функция, в окрестности замкнутом  круга радиуса $R$
с центром в нуле. В предшествующей первой части нашей работы были получены общие оценки на интеграл от положительной части радиальной максимальной характеристики роста  ${\mathsf M}_U(t):=\sup\bigl\{U(z)\bigm| |z|=r\bigr\}$  по возрастающей функции интегрирования $m$ на отрезке $[0,r]$ через разностную характеристику Неванлинны и модуль непрерывности функции $m$. Вторая часть работы даёт явный вид для  
таких оценок при условии, что модуль непрерывности функции $m$ не превышает  некоторую дифференцируемую функцию $h$ на открытом интервале $(0,r)$ с  единственным условием конечности точной верхней грани  $\sup\limits_{t\in (0,r)}\dfrac{h(t)}{th'(t)}<+\infty$. Этому условию удовлетворяют  любые степенные функции $t\mapsto t^d$ степени $d>0$. Оценки в определённом смысле оптимальны.

Библиография:  4 названия 

Ключевые слова: $\delta$-субгармоническая функция, радиальная макси\-м\-а\-льная характеристика роста, разностная характеристика Неванлинны, заряд Рисса, модуль непрерывности, условие Дини 
\end{abstract}

\markright{Мероморфные функции и разности субгармонических функций  \dots. II}


\section{Введение}

\subsection{Один недавний предшествующий результат}\label{s10}
Сохраняем и, по возможности, повторяем  определения и обозначения предшествующей первой части \cite{Kha21I} этой работы. Прежде приведём результат из  \cite[основная теорема]{Kha20}.

По-прежнему,  $\mes$ --- {\it линейная мера Лебега\/} на  вещественной оси $\RR$, а подмножество в $\RR$ или расширенную числовую функцию на подмножестве из $\RR$ со значениями  в расширенной вещественной оси $\overline \RR:=\RR\cup \pm\infty$ называем измеримыми, если они $\mes$-измеримы,  а  $\mes E:=\mes (E)$ для измеримого $E\subset \RR$,     
$$D(r):=\bigl\{z \in \CC \bigm| |z|< r\bigr\}, \quad
\overline  D(r):=\bigl\{z \in \CC \bigm| |z|\leq r\bigr\}, \quad 
\partial \overline  D(r):=\overline  D(r)
\setminus  D(r)
$$ 
--- соответственно {\it открытый\/} и  {\it замкнутый круги,\/} а также  {\it  окружность\/}  в комплексной плоскости   $\CC$ {\it радиуса $r\in \overline \RR^+$ с центром в нуле.\/}

Для $R\in \overline \RR^+:=\bigl\{x\in \overline \RR\bigm| x\geq 0\bigr\}$ и функции $v\colon D(R)\to \overline \RR$ 
  \begin{equation*}
{\sf M}_v(r):=\sup\bigl\{v(re^{i\theta})\bigm| \theta \in [0,2\pi)\bigr\} , \quad 0\leq r<R,
\end{equation*}
--- {\it максимальная характеристика роста функции $v$  на окружностях  $\partial \overline D(r)$,\/}
\begin{equation*}\label{{MC}C}
{\sf C}_v(r):=\frac{1}{2\pi}\int_0^{2\pi} v(re^{i\varphi})\dd \varphi 
\end{equation*}
---  {\it среднее по окружности\/ $\partial \overline{D}(r)$ функции $v$.}

 Для $ 0\leq r\leq R\in \RR^+ $ и   меры Бореля  $\mu$ на $\overline D(R)$ возрастающая функция 
\begin{subequations}\label{murad}
\begin{align}
\mu^{\rad}(r)&:=\mu^{\rad}\bigl(\overline D(r)\bigr)
\text{ при $r\in [0, R]$}
\tag{\ref{murad}$\mu$}\label{{murad}m}
 \\
\intertext{--- {\it  радиальная считающая функция меры $\mu$}, а } 
{\sf N}_{\mu}(r,R)&\overset{\eqref{{murad}m}}{:=}\int_{r}^{R}\frac{\mu^{\rad}(t)}{t}\dd t\in \overline \RR^+
\tag{\ref{murad}N}\label{{murad}N}
\end{align}
\end{subequations} 
--- {\it  разностная  усреднённая,\/} или  {\it проинтегрированная,  радиальная считающая функция меры $\mu$}
от двух переменных $0\leq r<R\leq +\infty$. 

Пусть  $U=u-v$ --- разность пары  субгармонических функций 
 $u\not\equiv -\infty$, $v\not\equiv -\infty$ в окрестности круга $\overline D(R)$ с {\it мерами Рисса\/} 
соответственно $\varDelta_u$ и  $\varDelta_v$. Таким образом,   $U\not\equiv\pm\infty$ ---  {\it $\delta$-суб\-га\-р\-м\-о\-н\-и\-ч\-е\-с\-к\-ая 
 функция\/} с  {\it зарядом  Рисса\/}  $\varDelta_U:=\varDelta_u-\varDelta_v$.
Различные  эквивалентные формы определения таких функций и их основные свойства приводятся и исследуются в \cite{Arsove53}, \cite[3.1]{KhaRoz18}. {\it Разностная   характеристика Неванлинны\/} такой   функции $U$ использовалась в нашей статье \cite{Kha20} и может быть определена 
как  функция двух переменных 
\begin{equation}\label{T}
{\sf T}_U(r,R)={\sf C}_{U^+}(R)-{\sf C}_{U^+}(r)+
{\sf N}_{\varDelta_U^-}(r,R),
\quad 0<r< R\in \RR^+, 
\end{equation}
где $U^+:=\sup\{0, U\}$ --- положительная часть функции $U$, а положительная мера $\varDelta_U^-:=\sup\{\varDelta_v, \varDelta_u\}-\varDelta_u\geq 0$ --- это {\it нижняя вариация\/} заряда Рисса  $\varDelta_{U}=\varDelta_u-\varDelta_v$ функции $U$. 

Для измеримых  $E\subset \RR$ и   $g\colon E\to \overline \RR$ наряду с существенной верхней гранью  
\begin{equation}\label{Lin}
\|g\|_\infty:=\inf \Bigl\{a\in \RR\Bigm| \mes\Bigl(\bigl\{x\in E\bigm|
g(x)>a\bigr\}\Bigr)=0 \Bigr\}
\text{ на $E$},
\end{equation}

используем и $L^p$-полунорму  функции $g$ на $E$  
\begin{equation}\label{Lpn}
 \|g\|_p:=\sqrt[p]{\int_E |g|^p\dd\, \mes}\quad \text{при  $1\leq p\in \RR^+$ на $E$}.
\end{equation}

\begin{ThSIw}[{\cite[основная теорема]{Kha20}}]\label{thm} 
Если $0< r_0< r<+\infty$, $1<k\in \RR^+$,  $E\subset [0,r]$ и $g\colon E\to \overline \RR$ измеримы, 
 $1<p\leq \infty$,  $1/p+1/q=1$,  $U\not\equiv \pm\infty$  --- $\delta$-субгармоническая функция на  $\CC$, 
а  $u\not\equiv -\infty$ --- субгармоническая функция на $\CC$, то  
\begin{subequations}\label{1m}
\begin{align}
\int_{E} {\sf M}_{U}^+(t)g(t)\!\dd t&\leq
\frac{4qk}{k-1} \bigl({\sf T}_{U}(r_0,kr)+{\sf C}_{U^+}(r_0)\bigr) \|g\|_p
\sqrt[q]{\mes E} \ln\frac{4kr}{\mes E}, \tag{\ref{1m}T}\label{inDl+}
\\ 
\int_{E} {\sf M}_{|u|}(t)g(t)\!\dd t&\leq
\frac{5qk}{k-1} \bigl({\sf M}_{u^+}(kr)+{\sf C}_{(-u)^+}(r_0)\bigr) \|g\|_p
\sqrt[q]{\mes E}\ln\frac{4kr}{\mes E}.
\tag{\ref{1m}M}\label{uM}
\end{align}
\end{subequations}
\end{ThSIw}

Нетрудно видеть, что  теорема о малых интервалах с весом при $p:=\infty$ содержит в себе и уточняет  все сформулированные во введении  к \cite[1.1]{Kha21I} предшествующие результаты. Из основного гораздо более общего результата настоящей второй части нашей работы  она будет выведена
в п. \ref{Ss2_5} следующего  \S~\ref{S12}, где этот общий результат и сформулирован. 

В связи с видом круглой скобки в правой части \eqref{inDl+}  с двумя слагаемыми далее будет удобнее  использовать именно такую форму {\it разностной характеристики Невалинны,\/} которую  можно определить через  ${\sf T}_{U}$ из \eqref{T} в виде 
\begin{subequations}\label{rT}
\begin{align}
{\boldsymbol T}_U(r,R)&:={\sf T}_{U}(r,R)+{\sf C}_{U^+}(r)
\overset{\eqref{T}}{=}{\sf C}_{U^+}(R)+{\sf N}_{\varDelta_U^-}(r,R),
\tag{\ref{rT}T}\label{{rT}T}
\\
\intertext{при любых $0<r<R\in \RR^+$, где правая часть позволяет определить и}
{\boldsymbol  T}_U(R)&:={\boldsymbol  T}_U(0,R)\overset{\eqref{{murad}N}}{:=}{\sf C}_{U^+}(R)+{\sf N}_{\varDelta_U^-}(0,R)\in \overline \RR^+.
\notag 
\end{align}
\end{subequations}

\subsection{Формулировки основного  результата из первой части работы}\label{S2}

Потребуется основное в \cite{Kha21I}
\begin{definition}\label{DefhR}
{\it Возрастающей функции\/}   $m\colon [0,r]\to  \RR$ {\it полной вариации\/}
\begin{equation}\label{{hR}wm}
 {\tt M} :=m(r)-m(0) \in \RR^+
\end{equation}
с  {\it модулем непрерывности\/} 
\begin{equation}\label{{hR}h}
\omega_m(t)\underset{t\in \RR^+}{:=}\sup\bigl\{ m(x)-m(x')\bigm|x-x'\leq t, \, 0\leq x'\leq x\leq r  \bigr\}\overset{\eqref{{hR}wm}}{\subset} [0,{\tt M}]
\end{equation}
сопоставляется  {\it диаметр стабилизации\/}
\begin{equation}\label{{hR}R}
{\sf d}_m:=\inf\bigl\{t\in \RR^+\bigm| {\omega}_{m}(t)= {\tt M}\bigr\}=\inf {\omega}_{m}^{-1}({\tt M})\leq r.
\end{equation}
\end{definition}

\begin{maintheorem}[{(\rm \cite[основная теорема]{Kha21I})}]\label{th1} Если\/  
$m\colon [0,r]\to\RR$ --- возрастающая функция с  модулем непрерывности, удовлетворяющим  условию Дини
\begin{equation}\label{{hR}i}
\int_0^{4R}\frac{{\omega}_{m}(t)}{t}\dd t<+\infty,
\end{equation}
то  для любой $\delta$-субгармонической функции $U\not\equiv\pm\infty$ на круге 
$\overline D(R)$ радиуса $R>r$ существует интеграл Лебега\,--\,Стилтьеса   с верхней оценкой 
\begin{equation}\label{U}
\int_0^r {\sf M}_{U}^+(t)\dd m(t)\leq \frac{6R}{R-r}
{\boldsymbol   T}_U(r,R) \max\biggl\{{\tt M}, \int_0^{{\sf d}_m}\ln \frac{4R}{t}\dd {\omega}_{m}(t)\biggr\}, 
\end{equation}
 где первый аргумент $r$ в ${\boldsymbol   T}_U(r,R)$ 
можно заменить на любое  $r_0\in [0,r]$, а   последний  интеграл Римана\,--\,Стилтьеса в \eqref{U} под операцией $\max$ --- на сумму
\begin{equation}\label{kint}
\int_0^{{\sf d}_m} \frac{{\omega}_{m}(t)}{t}\dd t+
{\tt M}\ln \frac{4R}{{\sf d}_m}\geq \int_0^{{\sf d}_m} \ln\frac{4R}{t}\dd {\omega}_{m}(t).
\end{equation} 
\end{maintheorem}

\section{Основной результат с явными оценками и примером}\label{S12} 

\subsection{Формулировка основного результата}
Если модуль непрерывности ${\omega}_{m}$ мажорируется некоторой дифференцируемой функцией из широкого класса, 
включающего в себя, в частности, все степенные функции строго положительной степени, то верхнюю оценку  из  \eqref{U} можно заменить на явную.
\begin{theorem}\label{th2}
Пусть $0<r\in \RR^+$, $h\colon [0,r]\to \RR^+$  --- непрерывная функция с  $h(0)=0$, дифференцируемая на $(0,r)$ и    удовлетворяющая условию   
\begin{equation}\label{{chrh}C}
{\sf s}_h:=\sup_{0< t< r} \frac{h(t)}{th'(t)}<+\infty.
\end{equation}
Тогда $h$ строго возрастает,  а  для любой возрастающей  функции $m\colon [0,r]\to\RR$ с полной вариацией  
${\tt M}\overset{\eqref{{hR}wm}}{:=}m(r)-m(0) \in \RR^+$ и 
модулем непрерывности 
\begin{equation}\label{{chrh}h}
{\omega}_m(t)\overset{\eqref{{hR}h}}{\leq} h(t)\text{ при всех $t\in [0,r]$}
\end{equation}
существует единственный прообраз  $h^{-1}({\tt M})\leq r$,  с которым   для  любой   $\delta$-суб\-г\-а\-р\-м\-о\-н\-и\-ч\-е\-с\-к\-ой  функции   $U\not\equiv \pm \infty$ в окрестности  замкнутого круга  $\overline D(R)$ радиуса   $R>r$ существует  интеграл Лебега\,--\,Стилтьеса с верхней оценкой   
\begin{equation}\label{Uh}
\int_0^r {\sf M}_{U}^+(t)\dd m(t) \leq \frac{6R}{R-r}{\boldsymbol  T}_U( r, R)\,
 {\tt M}\ln \frac{4e^{{\sf s}_h}R}{h^{-1}({\tt M})}, 
\end{equation}
где первый аргумент $r$ в ${\boldsymbol   T}_U(r,R)$ 
можно заменить на любое число $r_0\in [0,r]$. 
\end{theorem}

\subsection{О виде оценки в теореме \ref{th2}}
Определённую  оптимальность правой части оценки  \eqref{Uh} именно с логарифмическим множителем иллюстрирует
\begin{example}
Пусть $r:=2<R:=4$. Для каждого $s\in (0,1)$ рассмотрим возрастающую непрерывную функцию 
\begin{equation}\label{ms}
m(t)=\begin{cases}
1-s&\text{при $0\leq t\leq 1-s$},\\
t&\text{при $1-s\leq t\leq 1+s $},\\
1+s&\text{при $1+s\leq t\leq 2$}.
\end{cases}
\end{equation}
полной вариации ${\tt M}=2s$ с модулем непрерывности 
\begin{equation*}
\omega_{m}(t)=\begin{cases}
t&\text{при $0\leq t\leq 2s$}\\
2s&\text{при $t\geq 2s$}
\end{cases}, \qquad \omega_{m}(t)\leq t=:h(t), \quad h^{-1}(x)=x,  
\end{equation*}
и с ${\sf s}_h\overset{\eqref{{chrh}C}}{=}1$, что в данном случае даёт
\begin{equation}\label{tM4}
{\tt M}\ln \frac{4e^{{\sf s}_h}R}{h^{-1}({\tt M})}={\tt M}\ln \frac{4eR}{{\tt M}}=2s\ln \frac{8e}{s}.
\end{equation}
Рассмотрим мероморфную функцию $z\underset{z\in \CC}{\longmapsto} \dfrac{5}{z-1}$  и соответствующую ей  супергармоническую на $\CC$ и положительную на $\overline D (4)$ функцию 
\begin{equation}\label{U5}
U(z)\underset{z\in \CC}{=}\ln \frac{5}{|z-1|}, \quad\text{для которой   ${\sf M}_{U^+}(t)=\ln \frac{5}{|t-1|}$ при $t\leq 4$,}
\end{equation}
откуда для функции $m$ из \eqref{ms} получаем 
\begin{equation}\label{itM}
\int_0^2M_{U^+}(t)\dd m(t)=\int_{1-s}^{1+s}\ln \frac{5}{|t-1|}\dd t
=\int_0^s\ln\frac{5}{x}\dd x=s\ln\frac{5e}{s}.
\end{equation}
При этом для супергармонической функции $U$ из \eqref{U5} заряд Рисса $\varDelta_U$ --- это отрицательная масса $-1$ с носителем в точке $1\in \RR$, откуда    
\begin{equation*}
{\sf N}_{\varDelta_U^-}(0,4)=\int_1^4\frac{1}{t}\dd t=\ln 4, \quad 
 {\sf C}_{U^+}(4)=\frac{1}{2\pi}\int_0^{2\pi}
\ln \frac{5}{|4e^{i\theta}-1|}\dd \theta=\ln\frac{5}{4}
\end{equation*}
и, стало быть, ${\boldsymbol T}_U(4)={\sf C}_{U^+}(4)+{\sf N}_U(0,4) =\ln 5$. 
Отсюда правая часть неравенства   \eqref{Uh} при выборе $r_0:=0$  равна 
\begin{equation}\label{rs}
\frac{6R}{R-r}{\boldsymbol  T}_U( R)\,
 {\tt M}\ln \frac{4e^{{\sf s}_h}R}{h^{-1}({\tt M})}\overset{\eqref{tM4}}{=}
 \frac{6\cdot 4}{4-2} \ln 5 \cdot 2s\ln \frac{8e}{s}=24\ln5 \cdot s\ln\frac{8e}{s} 
\end{equation}
а интеграл из левой  части неравенства   \eqref{Uh} равен
\begin{equation}\label{ls}
\int_0^r {\sf M}_{U}^+(t)\dd m(t)\overset{\eqref{itM}}{=}s\ln\frac{5e}{s}.
\end{equation}
Сравнение правых частей \eqref{rs} и \eqref{ls} с варьированием $s\in (0,1)$ показывает, 
что неравенство   \eqref{Uh}  {\it должно содержать логарифмическую добавку\/} справа и в некотором смысле оптимально  по вкладу от $m$ с точностью до постоянных.  
\end{example}  


\subsection{Вывод теоремы о малых интервалах с весом из введения}\label{Ss2_5}
Установим  лишь  \eqref{inDl+}, поскольку  \eqref{uM} легко следует из него \cite[(24)--(26)]{Kha20}.
Переходя от функции $g$ к её положительной части $g^+$, при доказательстве можно, не умаляя общности, считать функцию $g$ положительной на $E\subset [0,r]$, равной нулю на $\RR\setminus E$ и с $\|g\|_p\neq 0$. 
Рассмотрим возрастающую функцию 
\begin{equation}\label{mg}
m(t):=\int_0^t g(s) \dd s=\int_0^t g(s)\mathbf{1}_E(s) \dd s, \quad  t\in [0,r],
\end{equation}
где $\mathbf{1}_E(t):=\begin{cases}
1 \quad\text{при $t\in E$},\\
0 \quad\text{при $t\notin E$},
\end{cases}$ --- характеристическая функция множества   $E$. 

Используя  неравенство Гёльдера, оцениваем её полную вариацию 
\begin{equation}\label{MM}
{\tt M}=\int_0^r g(s) \mathbf{1}_E(s) \dd s\leq \|g\|_p\biggl(\int_0^r\mathbf{1}_E^q(s)\dd s\biggr)^{1/q}
=  \|g\|_p\sqrt[q]{\mes E}\leq \|g\|_p\sqrt[q]{r} 
\end{equation}
 и  модуль непрерывности 
\begin{multline*}
\omega_m(t)\leq \sup_{x\in \RR}\int_{x}^{x+t} g(s)\mathbf{1}_E(s)\dd s\leq 
\|g\|_p  \sup_{x\in \RR}\biggl(\int_x^{x+t}\mathbf{1}_E^q(s)\dd s\biggr)^{1/q}
\\
\leq \|g\|_p  \sup_{x\in \RR}\biggl(\int_x^{x+t}\dd s\biggr)^{1/q}
=\|g\|_p t^{1/q}=:h(t) \quad \text{при всех $t\in \RR^+$}.
\end{multline*} 
Такая  функция $h$ удовлетворяет условиям теоремы \ref{th2} с 
\begin{equation}\label{h-1+}
{\sf s}_h\overset{\eqref{{chrh}C}}{=}q, \quad h^{-1}({\tt M})=\Bigl(\frac{\tt M}{\|g\|_p}\Bigr)^q.
\end{equation}
Отсюда по теореме \ref{th2} для $\delta$-суб\-г\-а\-р\-м\-о\-н\-и\-ч\-е\-с\-к\-ой  функции   $U\not\equiv \pm \infty$ в окрестности круга   $\overline D(R)$ радиуса   $R:=kr>r$ получаем 
\begin{multline*}
\int_0^r {\sf M}_{U}^+(t)g(t)\dd t\overset{\eqref{mg}}{=}\int_0^r {\sf M}_{U}^+(t)\dd m(t) \leq \frac{6k}{k-1}{\boldsymbol  T}_U( r, kr)\,
 {\tt M}\ln \frac{4e^{{\sf s}_h}R}{h^{-1}({\tt M})}\\
\overset{\eqref{h-1+}}{=}\frac{6k}{k-1}{\boldsymbol  T}_U( r, kr)\,
 {\tt M}\ln \frac{4e^{q}kr}{\bigl({\tt M}/{\|g\|_p}\bigr)^q}=
\frac{6qk}{k-1}{\boldsymbol  T}_U( r, kr)\,
 {\tt M}\ln \frac{e{\|g\|_p}\sqrt[q]{r}(4k)^{1/q}}{{\tt M}}
\end{multline*}
Для $0<b\in \RR^+$ функция $x\underset{x\in \RR^+}{\longmapsto} x\ln \frac{eb}{x}$ возрастающая на  $[0,b]$, и по  
\eqref{MM} оба вхождения ${\tt M}$ в правую часть можем заменить  на $\|g\|_p\sqrt[q]{\mes E}$, что даёт
\begin{multline*}
\int_0^r {\sf M}_{U}^+(t)g(t)\dd t \leq \frac{6qk}{k-1}{\boldsymbol  T}_U( r, kr)\,
\|g\|_p\sqrt[q]{\mes E}\ln \frac{e{\|g\|_p}\sqrt[q]{r}(4k)^{1/q}}{\|g\|_p\sqrt[q]{\mes E}}
\\=\frac{6k}{k-1}{\boldsymbol  T}_U( r, kr)\,
\|g\|_p\sqrt[q]{\mes E}\ln \frac{e^q4kr}{\mes E}\leq 
\frac{6qk}{k-1}{\boldsymbol  T}_U( r, kr)\,
\|g\|_p\sqrt[q]{\mes E}\ln \frac{4ekr}{\mes E}\\
\overset{\eqref{{rT}T}}{\leq} \frac{6qk}{k-1} \bigl({\sf T}_{U}(r_0,kr)+{\sf C}_{U^+}(r_0)\bigr) \|g\|_p
\sqrt[q]{\mes E} \ln\frac{4ekr}{\mes E} \quad\text{при всех $r_0\in (0,r]$}.
\end{multline*}
где учтено, что  $q\geq 1$. Таким образом, получена оценка  \eqref{inDl+}, правда, с увеличением   множителя $4$ до $6$ и появлением множителя $e$ в последней дроби, что  в данной тематике совершенно несущественно и вызвано лишь гораздо более общим характером теоремы \ref{th2}.

\section{Доказательство теоремы 1} Из условия \eqref{{chrh}C} сразу следует, что $h'>0$ на $(0,r)$, откуда  непрерывная  на $[0,r]$ функция $h$ {\it строго возрастает на\/} $[0,r]$. 
Очевидно, обратная функция $h^{-1}$ также {\it строго возрастает\/} и из условия \eqref{{chrh}h} следует
$h^{-1}\bigl({\omega}_m(t)\bigr)\leq t$ при всех $t\in [0, r]$. 
 Отсюда при  любом $t>{\sf d}_m$ по определению диаметра стабилизации ${\sf d}_m$ 
в  \eqref{{hR}R} получаем  ${\omega}_m(t)={\tt M}$ и, соответственно,   $h^{-1}({\tt M})\leq t$. 
В силу произвола в выборе $t>{\sf d}_m$ тогда
 \begin{equation}\label{hrr}
h^{-1}({\tt M})\leq {\sf d}_m\overset{\eqref{{hR}R}}{\leq} r.
\end{equation}
 Из условия  \eqref{{chrh}C} имеем  
\begin{equation}\label{ih'}
\int_0^{x}\frac{h(t)}{t}\dd t
\leq {\sf s}_h\int_0^x h'(t)\dd t
={\sf s}_h h(x) <+\infty\quad\text{для каждого $x\in [0,r]$}.
\end{equation}
Отсюда по условию \eqref{{chrh}h} следует условие \eqref{{hR}i} основной теоремы, поскольку
\begin{equation*}
\int_0^{4R}\frac{{\omega}_m(t)}{t}\dd t\overset{\eqref{ih'}}{\leq} \int_0^{r}\frac{h(t)}{t}\dd t
+\int_r^{4R}\frac{{\tt M}}{t}\dd t\leq {\sf s}_h h(r)+{\tt M}\ln \frac{4R}{r}<+\infty.
\end{equation*} 
Значит выполнено неравенство \eqref{U}, в котором правая часть \eqref{U} до ${\tt M}$
совпадает с фрагментом правой части  \eqref{Uh} до ${\tt M}$, а последний интеграл из \eqref{U} можно заменить на сумму из \eqref{kint} вида 
\begin{multline*}\label{kint=-}
\int_0^{{\sf d}_m} \frac{{\omega}_{m}(t)}{t}\dd t+
{\tt M}\ln \frac{4R}{{\sf d}_m}\overset{\eqref
{hrr}}{=}
 \biggl(\int_0^{h^{-1}({\tt M})}+\int_{h^{-1}({\tt M})}^{{\sf d}_m}\biggr) \frac{{\omega}_{m}(t)}{t}\dd t+
{\tt M}\ln \frac{4R}{{\sf d}_m}\\
\overset{\eqref{{chrh}h}}{\leq} 
\int_0^{h^{-1}({\tt M})}\frac{h(t)}{t}\dd t
+{\omega}_{m}({\sf d}_m)\int_{h^{-1}({\tt M})}^{{\sf d}_m} \frac{1}{t}\dd t+
{\tt M}\ln \frac{4R}{{\sf d}_m}
\\
\overset{\eqref{{hR}h}}{\leq} 
\int_0^{h^{-1}({\tt M})}\frac{h(t)}{t}\dd t+{\tt M}\ln \frac{{{\sf d}_m}}{h^{-1}({\tt M})}
+{\tt M}\ln \frac{4R}{{{\sf d}_m}}
\\\overset{\eqref{ih'}}{\leq}
{\sf s}_h h\bigl(h^{-1}({\tt M})\bigr)+{\tt M}\ln \frac{4R}{h^{-1}({\tt M})}
={\sf s}_h {\tt M}+{\tt M}\ln \frac{4R}{h^{-1}({\tt M})}=
{\tt M}\ln \frac{4e^{{\sf s}_h}R}{h^{-1}({\tt M})}.
\end{multline*} 
Таким образом,  для максимума  из \eqref{U} согласно  \eqref{kint} имеем 
\begin{equation*}
 \max\biggl\{{\tt M}, \int_0^{{\sf d}_m}\ln \frac{4R}{t}\dd {\omega}_{m}(t)\biggr\}
\leq  \max\biggl\{{\tt M},{\tt M}\ln \frac{4e^{{\sf s}_h}R}{h^{-1}({\tt M})}\biggr\} 
\overset{\eqref{hrr}}{\leq} {\tt M}\ln \frac{4e^{{\sf s}_h}R}{h^{-1}({\tt M})},
\end{equation*}
поскольку $R>r\overset{\eqref{hrr}}{\geq}h^{-1}({\tt M})$, что  вместе с общей оценкой \eqref{U} даёт  неравенство \eqref{Uh}
и завершает доказательство теоремы \ref{th2}.

\end{fulltext}

\end{document}